\documentclass{amsart}
\usepackage{amssymb,fancybox,natbib,amsmath}
\newcommand{\comment}[1]{\fbox{{\tiny\fbox{Skip:} }#1}}
\newcommand{\Lcomment}[1]{\fbox{\begin{minipage}{5in}{\tiny\fbox{Skip this:}}\\#1\end{minipage}}}
\renewcommand{\comment}[1]{}
\renewcommand{\Lcomment}[1]{}
      \newtheorem{theorem}{Theorem}[section]
       \newtheorem{proposition}[theorem]{Proposition}
       
       \newtheorem{lemma}[theorem]{Lemma}
\theoremstyle{remark}
       \newtheorem{remark}[theorem]{Remark}
         
\theoremstyle{definition}

\newtheorem{assumption}[theorem]{Assumption}

\newtheorem{note}{}

\def\E{{\mathbb E}}
\def\V{{\rm Var}}
\newcommand{\RR}{\mathbb{R}}
\newcommand{\CC}{\mathbb{C}}

\newcommand{\QQ}{\mathbb{Q}}
\newcommand{\calF}{\mathcal{F}}
\newcommand{\calA}{\mathcal{A}}

\numberwithin{equation}{section}

\author{W\l odek Bryc}
\address{
Department of Mathematical Sciences\\
University of Cincinnati\\
PO Box 210025\\
Cincinnati, OH 45221--0025, USA\\
Wlodzimierz.Bryc@UC.edu
}

\title{Markov processes with  free-Meixner laws}


\keywords
{martingale, free Meixner laws, generator, classical version of free L\'evy process}
\begin{document}
\maketitle
\begin{abstract}
%The abstract should state briefly the purpose of the research, the principal results and major conclusions.
We study  a time-non-homogeneous Markov process which arose from free probability, and which also appeared in the study of stochastic processes with linear regressions and quadratic conditional variances.
Our main result is the explicit expression for the  generator  of the (non-homogeneous) transition operator acting on functions that extend analytically to complex domain.

The paper is self-contained and does not use free probability techniques.
\comment{Boxed text will disappear in final version}
\end{abstract}
%\end{frontmatter}
\section{Introduction}
%State the objectives of the work and provide an adequate background, avoiding a detailed literature survey or a summary of the results.

In this paper we study a special class of (non-homogeneous) Markov processes whose univariate law form a semigroup with respect to the so called free additive convolution of measures.  These processes arise as the "classical versions" of the corresponding non-commutative free-L\'evy processes in the sense  that their time-ordered moments coincide, see \citet[page 144]{Biane:1998}. The same class of Markov processes also appeared as one of the examples in the study of "quadratic harnesses", i.e. processes with linear regression and quadratic conditional variances under double-sided conditioning with respect to past and future. The paper  however is self-contained and does not rely on  free probability techniques or "quadratic harnesses", except for motivation or "inspiration". (For example, the expression for the martingale in Proposition \ref{P:martingale} came from papers of Biane and Anshelevich but in this paper we verify the martingale property by direct integration.) To avoid distracting the reader, motivation and connections with free probability and with "quadratic harnesses"  are discussed in a separate section  at the end of the paper.

The paper is organized as follows. In Section \ref{Sect:Results} we define the family of Markov processes and state our main results.
Section \ref{Sect:Ele} collects elementary integrals needed for the proofs.
The integrals are then used in the proofs of the main results in Section \ref{Sect:Proofs}.
%Section
In %appendix
\ref{Sect:note} we  discuss relations to previous results, including connections  to free probability.

\section{Results}\label{Sect:Results}
We consider a family of probability measures $\{P_{s,t}(x,dy): 0\leq s<t, x\in\RR\}$ on Borel sets of the real line which depend on two auxiliary parameters $\theta\in\RR$ and $\tau\geq 0$.
The definition is somewhat cumbersome due to the possible presence  of an atom which may occur at the points that are given parametrically as
\begin{equation}\label{y*}
a_*(t)=\begin{cases}
-t/\theta & \mbox{ if } \tau=0, \theta\ne 0,\\
-t\frac{\theta-\sqrt{\theta^2-4\tau}}{2\tau} & \mbox{ if } \tau>0, \theta>0,\\
-t\frac{\theta+\sqrt{\theta^2-4\tau}}{2\tau} & \mbox{ if } \tau>0, \theta<0.
\end{cases}
\end{equation}
Probability measures $P_{s,t}(x,dy)$ are specified by their absolutely continuous component and discrete components (there is no singular component).
The continuous component  is given by the density
\begin{equation}\label{P_st_cont}
\frac{1}{2\pi}\frac{(t-s)\sqrt{4(t+\tau)-(y-\theta)^2}}
{\tau(y-x)^2+\theta(t-s)(y-x)+tx^2+sy^2-(s+t)xy+(t-s)^2}, %\: 1_{(y-\theta)^2<4(t+\tau)},
\end{equation}
supported on $y$ from  the interval $[\theta-2\sqrt{t+\tau},\theta+2\sqrt{t+\tau}]$.
The discrete component of  $P_{s,t}(x,dy)$  is zero except for the following cases.
\begin{enumerate}
\item If $\tau=0$, $\theta\ne 0$, and $x=a_*(s)=-s/\theta$, then with  $b^+=\max\{b ,0\}$ the discrete part of $P_{s,t}(x,dy)$   is given by
$$ \frac{\left(1-t/\theta^2\right)^+}{1-s/\theta^2}\delta_{a_*(t)}.$$
In particular,  the discrete component is absent for $t\geq \theta^2$.
\item If $\tau>0$, $\theta^2>4\tau$ and $x=a_*(s)$, then the discrete part of $P_{s,t}(x,dy)$ is given by
$$
\frac{\left(1-\frac{t}{2\tau}\frac{|\theta|-\sqrt{\theta^2-4\tau}}{\sqrt{\theta^2-4\tau}}\right)^+}
{1-\frac{s}{2\tau}\frac{|\theta|-\sqrt{\theta^2-4\tau}}{\sqrt{\theta^2-4\tau}}}\delta_{a_*(t)}.
$$
In particular,  the discrete component is absent for $t\geq {2 \tau\sqrt{\theta^2-4\tau}}/({ |\theta|-\sqrt{\theta^2-4\tau}})$.
\end{enumerate}
The laws  $P_{0,t}(0,dy)$ are the free Meixner laws in Note A.\ref{Sect:free_Meixner}.

Family $\{P_{s,t}(x,dy): 0\leq s<t, x\in\RR\}$ forms transition probabilities of a Markov process. This fact is implicit in \cite{Biane:1998}, and explicit in \cite[Theorem 4.3]{Bryc-Wesolowski-03}. Here we give a different proof based on the integral transform in Lemma \ref{H-transform}.
\begin{proposition}\label{P:Markov}
For every $\theta\in\RR$ and $\tau\geq 0$, there exists a right-continuous with left limits (cadlag) Markov process  $(X_t:t\geq 0)$ with state space $\RR$, initial state   $X_0=0$, and such that  for $0\leq s<t$, $\Pr(X_t\in U|X_s)=P_{s,t}(X_s,U)$ with probability one.
\end{proposition}

  The univariate laws of $X_t$ are $P_{0,t}(0,dy)$; these are the  free-Meixner laws  in the title of the paper, see Note A.\ref{Sect:free_Meixner}.

Next we describe a class of martingales associated with Markov process $(X_t)$.   We introduce the natural filtration $\calF_t:=\sigma(X_s:s\leq t)$, $t\geq 0$.

\begin{proposition}\label{P:martingale} Fix $z\in\CC$ such that $\tau |z|^2<1$. If $(X_t:t\geq 0)$ is the Markov process introduced in Proposition \ref{P:Markov},  then the complex-valued process
\begin{equation}\label{exp_mart}
M_t=\frac{1}{1-z (X_t-\theta)+(t+\tau)z^2}
\end{equation}
is an $\calF_t$-martingale for $0\leq t< 1/|z|^2 - \tau$.
\end{proposition}
 It might be worth pointing out that  $(M_t)$  is not a martingale  for $t> 1/|z|^2 - \tau$, as then
 $$\E(M_t)=\frac{t+\tau}{(t+\tau)^2z^2+\theta z (t+\tau)+\tau}$$ depends on $t$. See also  Note A.\ref{S:Biane-Anshelevich}.

 To state our next result we need additional notation.
 By $w_{m,\sigma^2}$ we denote the Wigner's semicircle law of mean $m$ and variance $\sigma^2>0$, given by the  density
\begin{equation}\label{eq:free_gauss}
w_{m,\sigma^2}(dx)=\frac{\sqrt{4\sigma^2-(x-m)^2}}{2\pi \sigma^2}1_{|x-m|\leq 2\sigma}(x)dx.
\end{equation}
% We include the degenerated law $\delta_m$ as $w_{m,0}$.

For $t>0$,  we consider the "generator"
$$
L_t(f)(x)=\lim_{h\to 0^+} \int \frac{f(y)-f(x)}{h}P_{t,t+h}(x,dy),
$$
defined on bounded measurable functions $f$ such that the limit exists. Our goal is to derive the expression for $L_t(f)$ when  $f$ belongs to a certain family $\mathcal{A}_t$  which contains all functions that extend analytically to the entire complex plain $\CC$.

To define this family $\mathcal{A}_t$, we denote by $r_t$ the radius of the disk centered at $\theta$ that contains the support   of $X_t$. Depending of the values of parameter $t,\theta,\tau$, this radius is the larger of the expressions
$2\sqrt{t}$ or $|\theta+t/\theta|$ when $\tau=0$ or the larger of $2\sqrt{t+\tau}$ and $\left((t+2\tau)|\theta|+t\sqrt{\theta^2-4\tau}\right)/(2\tau)$ when $\tau>0$, see \eqref{y*}.
Then   $f\in \mathcal{A}_t$ if there is $\delta>0$ such that $z\mapsto f(z)$ is analytic in the disk $|z-\theta|<\frac{5}{4} (r_t+\delta)$.

We now state our main result.
\begin{theorem}\label{T1} Fix $t>0$. If $f\in \mathcal{A}_t$, then for $x\in \mbox{supp} (X_t)$,
\begin{equation}\label{generator}
(L_t f)(x)=\frac{\partial}{\partial x}\int_\RR \frac{f(y)-f(x)}{y-x} w_{\theta,t+\tau}(dy).
\end{equation}
%\comment{Do we know what happens outside of support?}

\comment{Expect the same answer for $\tau<0$, with $t+\tau>0$...}
\end{theorem}
We remark that \eqref{generator} can be viewed as an analog of "Ito's formula" for instantaneous functions:
if $f$ is analytic in $\CC$ then
$$f(X_t)-\int_0^t L_s(f)(X_s) ds$$
is a martingale with respect to $(\calF_t)$. We also remark that at an atom of $X_t$ one should take the derivative before evaluating \eqref{generator} at $x=a_*(t)$. Equivalently,
$$
(L_t f)(x)=\int_\RR \frac{f(y)-f(x)-(y-x)f'(x)}{(y-x)^2} w_{\theta,t+\tau}(dy).
$$

We do not know the generators for Markov processes that correspond to  more general free-L\'evy processes; we also do not know the generators for  the  $q$-Meixner processes in \cite{Bryc-Wesolowski-03} when $q\ne 0,\pm 1$.

\section{Elementary integrals and an auxiliary Markov process}\label{Sect:Ele}
For complex $a_1,a_2,a_3,a_4$ let
$$
\widetilde{f}(x;a_1,a_2,a_3,a_4)=\frac{\sqrt{1-x^2}}{(1+a_1^2-2a_1x)(1+a_2^2-2a_2x)(1+a_3^2-2a_3x)(1+a_4^2-2a_4x)}.
$$
%The following integration results can be derived by partial fractions decomposition followed by integration of $\sqrt{1-x^2}/(1+a^2-2ax)$; the latter integral can be computed by
%complex change of variable $x=z+1/z$ over the unit circle $|z|=1$ by residua.
\begin{lemma}\label{L3.1}
If $|a_1|,\dots,|a_4|<1$, then
\begin{equation}\label{AW_int}
\int_{-1}^1\:\widetilde{f}(x;a_1,a_2,a_3,a_4)\:dx=K(a_1,a_2,a_3,a_4)\;,
\end{equation}
where \begin{equation}\label{<1} K(a_1,a_2,a_3)=\frac{\pi}{2}(1-a_1a_2a_3a_4)
\prod_{1\leq i<j\leq 4}(1-a_ia_j)^{-1}\;.  \end{equation}
\end{lemma}
\begin{proof} This integral is known (see Note A.\ref{A7}), but assuming $a_1,\dots,a_4$ are all distinct we provide the main steps of evaluation for completeness. By partial fractions decomposition, we only need to integrate four expressions of the form
$$\frac{a_1^3}{\prod_{j=2}^4\left[(a_1-a_j)(1-a_1a_j)\right]}\frac{\sqrt{1-y^2}}{\left(1+a_1^2-2 a_1 y\right)}.$$
Substituting $y=\cos \alpha$ and using the fact that $|a|<1$ we get
\begin{multline*}
\int_{-1}^1\frac{\sqrt{1-y^2}}{1+a^2-2 a y} dy=\frac12\int_0^{2\pi}\frac{\sin ^2 \alpha}{(1-a e^{i\alpha})(1-a e^{-i\alpha})}
d\alpha\\
%=\frac{i}{8}\oint_{|z|=1}\frac{(z-1/z)^2}{(1-az)(z-a)}dz
=\frac{i}{8}\oint_{|z|=1}\frac{(z^2-1)^2}{(1-az)(z-a)z^2}dz =\frac{\pi}{4} \left(1-\frac{1}{a^2}\right)+ \frac{\pi}{4} \left(1+\frac{1}{a^2}\right)=\pi/2,
\end{multline*}
with the last integral evaluated by residua at $z=a$ and $z=0$. (The third singularity at $z=1/a$ is outside of the unit disk.) Summing the four expressions from the partial fractions decomposition we get \eqref{<1}.
\end{proof}
In general, the integral in \eqref{AW_int} diverges when the parameters are on the unit circle; but there are two exceptions that arise from cancellations with the roots of $\sqrt{1-x^2}$: one parameter can take one of the values $\pm 1$ or a pair $(a_i,a_j)$ of parameters  can take the value $(-1,1)$. In these two exceptional cases the integral is still given by \eqref{<1}  still holds, as can be seen by taking the limits.

The integral in \eqref{AW_int}  converges also if some of the parameters are outside of the unit disk.  Since $1+a^2-2ax= a^2(1+1/a^2-2 x/a)$, formula \eqref{AW_int} can be used to evaluate such an integral. For example,
 if
$|a_2|,|a_3|,|a_4|<1$ and $|a_1|>1$, then
\begin{equation}\label{AW-large}
\int_{-1}^1\:\widetilde{f}(x;a_1,a_2,a_3,a_4)\:dx= K(1/a_1, a_2,a_3,a_4)/a_1^2,
\end{equation}
with
\begin{equation}\label{>1}  \frac{ K(1/a_1, a_2,a_3,a_4)}{a_1^2}=\frac{\pi(a_1-a_2a_3a_4)}
{2(a_1-a_2)(a_1-a_3)(a_1-a_4)(1-a_2a_3)(1-a_2a_4)(1-a_3a_4)}\;.\end{equation}

\subsection{Probability measures}
We now introduce a two-parameter family of  probability measures with parameters that satisfy the following. \begin{assumption}\label{Asume1}
Let $a_1,a_2$ be either real or complex conjugate, such  that their product satisfies $a_1a_2<1$.
\end{assumption}
Assumption \ref{Asume1} is a concise way of stating that either $a_1=\bar a_2$ are from the unit disk of the complex plain, or $a_1,a_2$ are real and  at least one of them is in the interval $(-1,1)$,  or if both are real but outside of $(-1,1)$ then they have opposite signs. We will need to consider these cases separately in the definitions and in the proofs.

Under Assumption \ref{Asume1},  $\widetilde f(y;a_1,a_2,0,0)$ is real-valued, positive, and integrable.
To confirm this, we need to consider separately the case when $a_1=\bar a_2$, and the case when $a_1,a_2$ are real. To see positivity for real $a_1,a_2$,  we write
$$(1+a_1^2-2 a_1 y)(1+a_2^2-2 a_2 y)=\left|(1-a_1 e^{i\alpha_y})(1-a_2 e^{i\alpha_y})\right|^2$$
with $\alpha_y=\arccos y$.

The corresponding normalizing constant % \eqref{<1}
$$k(a_1,a_2)=K(a_1,a_2,0,0)=\frac{\pi}{2(1-a_1a_2)}$$
is well defined and positive.
We  therefore  introduce the non-negative  function
\begin{equation}\label{f-def}
f(y;a_1,a_2)=\frac{1}{k(a_1,a_2)}\widetilde f(y;a_1,a_2,0,0)1_{[-1,1]}(y).
\end{equation}
By \eqref{AW_int}, $f$ is a probability density function when $|a_1|,|a_2|<1$.  For other values of admissible parameters, it is easy to check that $f(y) dy$ is a sub-probability measure. Adding the missing mass as the weight of (carefully selected!) atoms, we
consider the following two-parameter family of probability measures:
\begin{equation}\label{f_density}
\nu(dy; a_1,a_2)=\begin{cases}
 \displaystyle  f(y;a_1,a_2)\,dy &\mbox{ if $|a_1|,|a_2|<1$},\\
\displaystyle \frac{(1\mp a_2)\sqrt{1\pm x}}{\pi\sqrt{1\mp x}(1+a_2^2-2a_2 x)} &\mbox{ if $-1<a_2<1, a_1=\pm 1$},\\
 \displaystyle \frac{1}{\pi\sqrt{1-x^2}} &\mbox{if $a_1=\pm 1$, $a_2=-a_1$},\\
 \displaystyle  f(y;a_1,a_2)\,dy +w(a_1,a_2)\delta_{y(a_1)} &\mbox{ if $-1<a_2<1$, $|a_1|>1$},\\
 \displaystyle  f(y;a_1,a_2)dy+w(a_1,a_2) \delta_{y(a_1)}+w(a_2,a_1)\delta_{y(a_2)}& \mbox{if $a_1>1$ and $a_2<-1$},
\end{cases}
\end{equation}
where the locations of the atoms are $y(a)=(a+1/a)/2$ and the weights of the atoms are
\begin{equation}\label{weight}
w(a,b)=\frac{a^2-1}{a^2-ab}.
\end{equation}
It is straightforward to verify that $0<w(a_1,a_2)<1$ and that
$$w(a_1,a_2)=1-\frac{k(1/a_1,a_2)}{a_1^2k(a_1,a_2)}= 1-\int_{-1}^1 f(x;a_1,a_2)dx$$
when $a_1,a_2$ are real, $a_1a_2<1$, $-1<a_2<1$ and $|a_1|>1$. Furthermore, it is clear  that  $w(a_1,a_2),w(a_2,a_1)>0$ and that
$$w(a_1,a_2)+w(a_2,a_1)=1+\frac{1}{a_1a_2}=1-\frac{k(1/a_1,1/a_2)}{a_1^2a_2^2k(a_1,a_2)}=1-\int_{-1}^1 f(x;a_1,a_2)dx$$
when $a_1>1,a_2<-1$.

We extend the definition \eqref{f_density} to the entire range of admissible parameters $a_1,a_2$ by symmetry: we
request that $\nu(dy;a_1,a_2)=  \nu(dy;a_2,a_1)$ also in all cases omitted from \eqref{f_density}.

We note the following elementary formulas.
\begin{proposition}\label{P-moms}
The mean of  $\nu(dy;a_1,a_2)$ is
$$m=\int_\RR y\,\nu(dy;a_1,a_2)=(a_1+a_2)/2,$$ and the variance  is
$$\int_\RR (y-m)^2\,\nu(dy;a_1,a_2)=(1-a_1a_2)/4.$$
 For $|z|<1$,
\begin{equation}\label{CS-transform}
\int_\RR \frac{1}{1+z^2 - 2 zy}\nu(dy;a_1,a_2)=\frac{1}{(1-a_1 z)(1-a_2 z)}.
\end{equation}
\end{proposition}
\begin{proof}
To compute the moments we take the derivatives of both sides of \eqref{CS-transform}  at $z=0$.

To derive formula \eqref{CS-transform} we need to consider separately each case that appears in \eqref{f_density}. In each case we apply \eqref{AW_int} to evaluate the integral over the absolutely continuous component of the measure, and add the corresponding contribution of the discrete component.

In the case $|a_1|,|a_2|<1$,  the left hand side of \eqref{CS-transform} is $K(a_1,a_2,z,0)/K(a_1,a_2,0,0)$. From \eqref{<1} we get \eqref{CS-transform}.

In the case $|a_1|>1$, $|a_2|<1$ we use \eqref{AW-large}. From the continuous part we get
 $$\frac{K(1/a_1,a_2,z,0)}{a_1^2K(a_1,a_2,0,0)}=\frac{1-a_1a_2}{(a_1-a_2)(1-a_2z)(a_1-z)}.
 $$
 The discrete part contributes
 $$\frac{w(a_1,a_2)}{1+z^2-2z y(a_1)}=
 \frac{w(a_1,a_2)}{(1-za_1)(1-z/a_1) }=
 \frac{a_1^2-1}{(a_1-a_2)(1-a_1 z)(a_1-z)}.$$
 The sum of these two contributions  gives the right hand side of \eqref{CS-transform}.

If $a_1>1$ and $a_2<-1$, the continuous part contributes
\begin{equation}\label{KKK}
\frac{K(1/a_1,1/a_2,z,0)}{a_1^2a_2^2K(a_1,a_2,0,0)}=-\frac{1}{(z-a_1)(z-a_2)}.
\end{equation}

 The discrete part  contributes
 $$\frac{w(a_1,a_2)}{1+z^2-2z y(a_1)}+\frac{w(a_2,a_1)}{1+z^2-2z y(a_2)}=
\frac{w(a_1,a_2)}{(1- a_1z)(1-z/a_1)}+\frac{w(a_2,a_1)}{(1-a_2z)(1-z/a_2)} $$
 $$
=\frac{\left(1+a_1 a_2\right) (z^2+1)-2 \left(a_1+a_2\right) z}{\left(z-a_1\right) \left(1-a_1z\right)
   \left(z-a_2\right) \left(1-a_2z\right)} =\frac{1}{\left(1-a_1z\right) \left(1-
   a_2z\right)}+\frac{1}{\left(z-a_1\right) \left(z-a_2\right)}.
 $$
The sum of this expression and \eqref{KKK}  gives  the right hand side of \eqref{CS-transform}.

The remaining cases with $a_1$ or $a_2$ taking values $\pm 1$ are the limits of the above.
\end{proof}
%\Lcomment{
%The $H$-transform for general free Askey-Wilson is
%$$\frac{a(1) a(2) a(3) a(4) z^2}{(z a(1)-1) (z a(2)-1) (z a(3)-1) (z a(4)-1)}
%   $$
%   $$+\frac{ \left(-a(2)
%   a(3) a(4) a(1)^2+\left(-a(3) a(4) a(2)^2-(a(3)+a(4)) (a(3)
%   a(4)-1) a(2)+a(3) a(4)\right) a(1)+a(2) a(3) a(4)\right) z}{(z a(1)-1) (z a(2)-1) (z a(3)-1) (z a(4)-1)
%   (a(1) a(2) a(3) a(4)-1)}
%   $$
%$$+\frac{1}{(z a(1)-1) (z a(2)-1) (z a(3)-1) (z a(4)-1)
%  }
%   $$

%
%   }

The following identity will be used to verify Chapman-Kolmogorov equations.
\begin{proposition}\label{P.m}
If $a_1,a_2$ satisfy Assumption \ref{Asume1} then for all $-1<m<1$, and all Borel sets $U$,
\begin{equation}\label{nu-m}
\nu(U;ma_1, ma_2)=\int_\RR \nu\left(U; m(x+\sqrt{x^2-1}),m(x-\sqrt{x^2-1})\right) \nu(dx; a_1,a_2).
\end{equation}
\end{proposition}
A short  proof uses the following $H$-transform.
\begin{lemma}\label{H-transform} A compactly supported probability measure $\nu$ is determined uniquely by the function
$z \mapsto H(z)=\int (1+z^2- 2 zy)^{-1} \nu(dy)$ for $z$ in a neighborhood of $0$.
\end{lemma}
\begin{proof} %$H(z)=G((z+1/z)/2)/(2z)$, where $G$ is the Cauchy transform.
A compactly supported measure is determined uniquely by its moments. The $k$-th moment of $\nu$ can be computed  from the $k$-th derivative of $H$ at $z=0$ and the moments of lower orders.
\end{proof}

\begin{proof}[Proof of Proposition \ref{P.m}]
Applying \eqref{CS-transform}  twice, the $H$-transform of the right hand side of \eqref{nu-m} is
$$
\int \frac{1}{(1-z m (x+\sqrt{x^2-1}))(1-z m (x-\sqrt{x^2-1}))}\nu(dx;a_1,a_2)
$$
$$
=\int \frac{1}{1+(m z)^2-2 (m z) x }\nu(dx;a_1,a_2)=\frac{1}{(1-mz a_1)(1-m za_2)}
.$$
From \eqref{CS-transform} we see that this matches the $H$-transform of the left hand side of \eqref{nu-m}.
\end{proof}

\subsection{An auxiliary Markov process}
Next we define transition probabilities of a Markov process with state space $\RR$ and time $T=(CD,\infty)$, where $C,D$ are either real  such that $CD\geq 0$  or complex conjugate.
%Without loss of generality we assume $|C|\leq |D|$. \comment{Do we use this?}

For $t\in(CD,\infty)$ we
define   probability measures
$$
\mu_t(dy)=\nu\left(dy;\frac{C}{\sqrt{t}},\frac{D}{\sqrt{t}}\right),
$$
and for $s<t$, $s,t\in[CD,\infty)$ and any real $x$ we define probability measures
$$
\mu_{s,t}(x,dy)=\nu\left(dy;\sqrt{\frac{s}{t}}(x+\sqrt{x^2-1}),\sqrt{\frac{s}{t}}(x-\sqrt{x^2-1})\right).
$$
Note that these measures are well defined: in each case the corresponding  parameters $a_1,a_2$ are either real or complex conjugates, and their product satisfies  $a_1a_2<1$.

We want to check that these measures form a Markov family, that is:
\begin{proposition}\label{P.3.2}
For $CD<s<t$,
\begin{equation}\label{marg}
\mu_t(dy)=\int_\RR\:\mu_{s,t}(x,dy)\mu_s(dx).
\end{equation}
For $CD<s<t<u$ and real $x$,
\begin{equation}\label{cond}
\mu_{s,u}(x,dz)=\int_\RR\:\mu_{t,u}(y,dz)\mu_{s,t}(x,dy).
\end{equation}
In addition, we have
\begin{equation}\label{pre-Martingale}
\int_\RR (1+z^2-2 z y)^{-1}\mu_{s,t}(x,dy)=\begin{cases}
\displaystyle\frac{t}{t+s z^2-2\sqrt{st}z x} &\mbox{if $|z|<1$},\\
\displaystyle\frac{t}{tz^2+s -2\sqrt{st}z x}&\mbox{if $|z|>1$}.
\end{cases}
\end{equation}

\end{proposition}
\begin{proof}
Formula \eqref{marg} follows from \eqref{nu-m} applied to $a_1=C/\sqrt{s}$, $a_2=D/\sqrt{s}$ and $m=\sqrt{s/t}$. Formula  \eqref{cond} follows from \eqref{nu-m} applied to $a_1=\sqrt{\frac{s}{t}}(x+\sqrt{x^2-1})$, $a_2=\sqrt{\frac{s}{t}}(x-\sqrt{x^2-1})$ and $m=\sqrt{t/u}$. Formula \eqref{pre-Martingale} follows from \eqref{CS-transform} applied  to $z$ when $|z|<1$ or to $1/z$ when $|z|>1$.
\end{proof}
\begin{remark}
The construction works also for real $C,D$ such that $CD<0$, with time   $T=(0,\infty)$.
\end{remark}

\section{Proofs of the main results}\label{Sect:Proofs}

%
%Lemma \ref{L3.1} becomes relevant once  we substitute
%\begin{equation}\label{eq:Y}
%Y_t:=2\sqrt{t+\tau}(X_{t+\tau}-\theta)
%\end{equation}
%and denote by $C,D$  the (possibly complex) roots of $1+\theta x+\tau x^2=0$.
%%\begin{equation} C+D=-\theta, \; CD=\tau.
%%\end{equation}

%The substitution transforms measures $P_{s,t}(x,dy)$  into measures
%$$
%p_{s,t}(x,dy)=\begin{cases}
%...
%\\
%f(y;0,0, \sqrt{\frac{s}{t}}(x-\sqrt{x^2-1}),\sqrt{\frac{s}{t}}(x+\sqrt{x^2-1}) &\mbox{ if $|x|<1$}
%\\
%\end{cases}
%$$

\begin{proof}[Proof of Proposition \ref{P:Markov}]
Let $C,D$ denote the roots of $z^2+\theta z+\tau=0$, so that  $\tau=CD$ and $\theta=-(C+D)$. Of course, $C,D$ are either real or complex conjugate,  so the Markov process  $(Y_t)_{t>\tau}$  from Proposition \ref{P.3.2} is well defined.

For rational $t>0$ define
\begin{equation}\label{eq:Y}
X_t=\theta+{2\sqrt{t+\tau}} {Y_{t+\tau}}.
\end{equation}
Then $(X_t)_{t\in\QQ_+}$ is a Markov process.
From Proposition \ref{P-moms} we see that
$$\E(Y_t)=\frac{C+D}{2\sqrt{t}}, \; \V(Y_t)=\frac{t-CD}{4t},$$
so $E(X_t)=0$ and $E(X_t^2)=t$.

From \eqref{pre-Martingale} with   $z$ replaced by $z\sqrt{t+\tau}$, we get
$$\E\Big(\frac{1}{1+z^2(t+\tau)-2z\sqrt{t+\tau}Y_{t+\tau}}\Big|Y_{s+\tau}\Big)=\frac{1}{1+z^2(s+\tau)-2z\sqrt{s+\tau}Y_{t+\tau}}$$
for all $s<t$  such that $t+\tau<1/|z|^2$.
This shows that
Proposition \ref{P:martingale} holds over positive rational $t$.
In particular, taking the derivative with respect to $z$ at $z=0$ we see that
 $\theta+M_t'(0)=X_t$ is a (square-integrable)  martingale. Therefore $X_t=\lim_{q\to t^+,q\in Q} X_q$ exists  almost surely, and defines a Markov process with right-continuous trajectories that have left limits, see \cite[Theorem 6.27]{Kallenberg:1997}. Of course, the transition probabilities of $(X_t)$ are re-calculated from the transition probabilities of $(Y_{t+\tau})$, and $X_0=0$ since $\V(X_t)=t$ for rational $t>0$.
(Details of calculation of transition probabilities for $(X_t)$ are omitted.)
  \end{proof}

\begin{proof}[Proof of Proposition \ref{P:martingale}]
We already saw that the result holds true for rational $t$. The general version follows by taking the limit.
\end{proof}
\begin{proof}[Proof of Theorem \ref{T1}]
Fix $f\in \calA_t$ such that $f$ is analytic in the disk $|u-\theta|<5/4 (r_{t}+2\delta)$ and take $h>0$ small so enough the support of $X_{t+h}$ is in the disk $|u-\theta|<r_t+\delta$. Let $\gamma$ be a curve in the first disk that encloses the support of $X_{t+h}$, and let
$x$ be in the support of $X_{t+h}$. Substituting $u=1/z+\theta+(t+\tau+h)z$ in the
Cauchy formula $f(x)=\frac{1}{2\pi i}\oint_{\gamma}{f(u)}(u-x)^{-1}du$, we get
\begin{equation}\label{g_t}
f(x)=\frac{1}{2\pi i}\oint_{|z|=1/(r_{t}+\delta)}
 \frac{g_{t+h}(z)}{1-z(x-\theta)+(t+h+\tau)z^2},
\end{equation}
where
$$
g_t(z)=\left((t+\tau) z-\frac1z\right)f\left(\theta+(t+\tau)z+\frac1z\right),
$$
and $\gamma$ is the ellipse $u(s)=\theta+ (r_t+\delta) e^{-i s} +\frac{t+h+\tau}{r_t+\delta}e^{i s}$.
Here we observe that
$$|u(s)-\theta|\leq r_t+\delta+ \frac{r^2_{t+h}}{4 (r_t+\delta)}<\frac{5}{4} (r_t+\delta) $$ for $h$ small enough, so $f$ is analytic in a disk that contains $\gamma$. Also $\gamma$ encloses the interval $(\theta-r_t-\delta,\theta+r_t+\delta)$ which for small enough $h\geq 0$ contains the support of  $X_{t+h}$.
Recall that $r_t\geq 2\sqrt{t+\tau}$. From \eqref{g_t} we see that by Proposition \ref{P:martingale}
applied with $h>0$ small enough so that $t+h+\tau<(r_t+\delta)^2$,
\begin{multline*}
L_t(f)(x)=\lim_{h\to 0^+} \frac{1}{2\pi i}\oint_{|z|=1/(r_{t}+\delta)}
 \frac{ (g_{t+h}(z)-g_t(z))/h}{1-z(x-\theta)+(t+\tau)z^2} dz
 \\=
 \frac{1}{2\pi i}\oint_{|z|=1/(r_{t}+\delta)}
 \frac{1}{1-z(x-\theta)+(t+\tau)z^2} \frac{\partial g_{t}(z)}{\partial t} dz.
\end{multline*}
Differentiating \eqref{g_t} with respect to $h$ at $h=0$ we get
\begin{multline*}
 \frac{1}{2\pi i}\oint_{|z|=1/(r_{t}+\delta)}
 \frac{1}{1-z(x-\theta)+(t+\tau)z^2} \frac{\partial g_{t}(z)}{\partial t} dz\\
 =
 \frac{1}{2\pi i}\oint_{|z|=1/(r_{t}+\delta)}
 \frac{z^2 g_{t}(z)}{(1-z(x-\theta)+(t+\tau)z^2)^2} dz.
\end{multline*}
So
\begin{equation}\label{LLL}
L_t(f)(x)=\frac{1}{2\pi i}\oint_{|z|=1/(r_{t}+\delta)}
 \frac{z^2 g_{t}(z)}{(1-z(x-\theta)+(t+\tau)z^2)^2} dz.
\end{equation}
We now verify that the right hand side of \eqref{generator} gives the same answer.
From \eqref{g_t} with $h=0$ we see that for $x,y$ in the support of $X_t$,
\begin{multline}\label{fraction}
\frac{f(y)-f(x)}{y-x}
=\frac{1}{2\pi i} \oint_{|z|=1/(r_{t}+\delta)}\frac{zg_t(z) dz}{(1-z(x-\theta)+(t+\tau)z^2)(1-z(y-\theta)+(t+\tau)z^2)}.
\end{multline}

Now we note that the support of the semicircle law $w_{\theta,t+\tau}$ is contained in the support of $X_t$, and that with $u=\sqrt{t+\tau}z$ in the unit circle,
by Proposition \ref{P:martingale} applied to the case of semicircle law, i.e., to $\theta=\tau=0$ we have
$$\int_\RR \frac{1}{1-z(y-\theta)+(t+\tau)z^2}w_{\theta,t+\tau}(dy)=
\int_\RR \frac{1}{1-uy+u^2}w_{0,1}(dy)=1.$$
Thus integrating \eqref{fraction} we get
$$\int_\RR \frac{f(y)-f(x)}{y-x} w_{\theta,t+\tau}(dy)=\frac{1}{2\pi i}
\oint_{|z|=1/(r_{t}+\delta)}\frac{zg_t(z)}{1-z(x-\theta)+(t+\tau)z^2}dz.$$
Taking the derivative of this expression with respect to $x$ and using  \eqref{LLL} we get \eqref{generator}.
\end{proof}

\subsection*{Acknowledgements}
We thank M. Anshelevich, M. Bo\.zejko,  %A. Hassairi,
P. Vallois, and J. Weso\l owski for helpful suggestions and discussions.  This research partially supported by NSF grant  DMS-0904720.

\bibliographystyle{apa}
\bibliography{free-gen}
\appendix
 \section{Notes}\label{Sect:note}
 \comment{These unpolished comments on relations to free probability, motivations, etc.}
%\subsection{}
\begin{note}\label{Sect:free_Meixner} Up to affine transformations, the univariate laws $\{P_{0,t}(0,dy):t>0\}$ of $(X_t)_{t>0}$  come from a two-parameter family of what is now called the "free Meixner laws".  These laws were introduced as the orthogonality measures of systems of polynomials with constant recursions in
%Anshelevich
 %Saitoh and Yoshida
\cite{Saitoh-Yoshida01} who found the explicit formula, analyzed free infinite divisibility and pointed out that this class includes a number of laws of interest in free probability; the term "free Meixner" was introduced in \cite{Anshelevich:2003}.  Further properties were studied in a series of papers
\cite{Anshelevich:2004,Anshelevich:2005,Anshelevich:2007,Anshelevich:2008,Bozejko-Bryc-04}.

The
free Meixner laws can be classified into six types: Wigner's semicircle (free Gaussian) which corresponds to our $\tau=\theta=0$, free
Poisson (also known as Marchenko-Pastur) which corresponds to our $\tau=0$, $\theta\ne 0$,  free Pascal (also known as free negative binomial) which corresponds to our $\theta^2>4\tau>0$, free Gamma which corresponds to our $\theta^2=4\tau>0$,  a law
that one may call pure free Meixner, and  the free binomial law which corresponds to the case $\tau<0$ that is not considered in this note; the complete list of cases builds on \cite[Remark 2.5 and Examples 3.4, 3.6]{Saitoh-Yoshida01}, \cite[Theorem 4]{Anshelevich:2003} and  appears in \cite[Theorem 3.2]{Bozejko-Bryc-04} or in \cite[Remark 4]{anshelevich2009bochner}.
\end{note}

\begin{note}\label{Sect:Markov}
The Markov process $(X_t)$   can be introduced as follows.
Except for the free binomial family, the free Meixner laws  are infinitely-divisible with respect to
the additive free convolution, see \cite[Theorem 3.2]{Saitoh-Yoshida01}, and are therefore univariate laws of  non-commutative free L\'evy processes. By   \citet[Theorem 3.1]{Biane:1998} there exists a unique Markov process $(X_t)$ with the same univariate laws and the same time-ordered joint moments (if they exist).
 In particular, if $\theta=\tau=0$, the univariate laws of $(X_t)$  are the semicircle laws $w_{0,t}$, $t>0$, and the corresponding transition probabilities   appear in  \cite[Section 5.3]{Biane:1998}.

The same family of Markov processes can be specified by conditional means and conditional variances, see \cite[Theorem 4.3]{Bryc-Wesolowski-03},
and our construction   is based on the formulas from that paper.

% Transition probabilities $P_{s,t}(x,dy)$  can be described more concisely but less intuitively by their
%Cauchy-Stieltjes transform
%%\begin{equation}
%\begin{eqnarray}
%\label{f-transition}
%&\displaystyle\int_{\RR}\frac{1}{z-y}P_{s,t}(x,dy)&  \\  \nonumber =&\displaystyle \frac12
%\frac{ (t+s+2\tau)(z - x)+(t-s)\theta  -(t-s)
%    \sqrt{ (z-\theta)^2 -
%        4 (t+\tau) }}{\tau (z-x)^2+\theta(t-s)(z-x)+tx^2+sz^2-(s+t)x z+(t-s)^2 }.&
%%\end{equation}
%\end{eqnarray}
\end{note} \begin{note}\label{S:Biane-Anshelevich}  Proposition \ref{P:martingale} can be deduced from \citet[Proposition 4.3.1]{Biane:1998}. However to do so when $\tau>0$ one needs to use a non-trivial substitution that appears in \cite[page 236]{Anshelevich:2003}. Additional analysis is needed to determine explicitly the allowed range of $t$ which  is  crucial for our proof of Theorem \ref{T1}.
%\subsection{}
% The Wigner's semicircle  law \eqref{eq:free_gauss} appears as the limiting spectrum of large random matrices. For free Meixner laws, it plays the role of L\'evy-Chinchine measure, see ...

\end{note}
\begin{note} For $\theta=\tau=0$,  formula  \eqref{generator} agrees with the non-commutative result \citet[page 392]{Biane:1998b} after correcting their expression by a factor of 2, and with \cite[page 150]{BKS:1997}, who consider a closely related classical Markov process
  $(e^t X_{e^{-2t}})_{t>0}$ with the generator
 $$
L_t(f)(x)=x f'(x)-2\frac{\partial}{\partial x} \int \frac{f(y)-f(x)}{y-x}w_{0,1}(dy).
$$
 \end{note}
 \begin{note}
 Generators of more general Markov processes that arise from free L\'evy processes can be read out from  \cite[Corollary 10]{anshelevich2002ito}. For properties of operator $f\mapsto \int \frac{f(y)-f(x)}{y-x}\mu(dy)$ with compactly supported $\mu$ see
 \cite[Proposition 1]{anshelevich2009bochner}.
\end{note}

 \begin{note} \label{Classical version}
Combining  Proposition \ref{P:martingale}  with  Lemma \ref{H-transform}, and \cite[Proposition 4.3.1]{Biane:1998} with \cite[page 236]{Anshelevich:2003} one verifies that the action of transition probabilities of $(X_t)$ on  polynomials $f$ coincides with the action of non-commutative conditional expectation, so the joint moments of our process $(X_t)$  indeed match the non-commutative moments as explained in \cite[page 161]{Biane:1998}
\end{note}

\begin{note}\label{A7} Lemma \ref{L3.1} is an elementary case of  the Askey-Wilson integral \cite[(2.1)]{Askey-Wilson-85}.    \citet[Eqn (1.3)]{Ismail:1995} state this elementary integral when $a_3=a_4=0$. \comment{Need a good reference!}
\end{note}
\begin{note}A version of Lemma \ref{H-transform} holds true also for non-compactly supported measures,
as $H(z)=G((z+1/z)/2)/(2z)$, where $G(u)=\int(u-x)\nu(dx)$ is the Cauchy-Stieltjes  transform of $\nu$.
\end{note}

\end{document}